# CORRECTION
# BROWNIAN MODELS OF OPEN PROCESSING NETWORKS: CANONICAL REPRESENTATION OF WORKLOAD

By J. Michael Harrison

Due to a printing error the above mentioned article had numerous equations appearing incorrectly in the print version of this paper. The entire article follows as it should have appeared. IMS apologizes to the author and the readers for this error.

A recent paper by Harrison and Van Mieghem explained in general mathematical terms how one forms an "equivalent workload formulation" of a Brownian network model. Denoting by $Z(t)$ the state vector of the original Brownian network, one has a lower dimensional state descriptor $W(t) = MZ(t)$ in the equivalent workload formulation, where $M$ can be chosen as any basis matrix for a particular linear space. This paper considers Brownian models for a very general class of open processing networks, and in that context develops a more extensive interpretation of the equivalent workload formulation, thus extending earlier work by Laws on alternate routing problems. A linear program called the static planning problem is introduced to articulate the notion of "heavy traffic" for a general open network, and the dual of that linear program is used to define a canonical choice of the basis matrix $M$. To be specific, rows of the canonical $M$ are alternative basic optimal solutions of the dual linear program. If the network data satisfy a natural monotonicity condition, the canonical matrix $M$ is shown to be nonnegative, and another natural condition is identified which ensures that $M$ admits a factorization related to the notion of resource pooling.

## Contents

---

**1. Introduction.** Brownian networks are a class of stochastic system models introduced in [6] and later used to approximate queueing networks of various kinds under conditions of heavy traffic [11, 12, 13, 14, 15, 16, 17, 18, 19]. To be more specific, Brownian networks arise as heavy traffic approximations of multiclass queueing networks in which system managers have a dynamic control capability. It has also been observed repeatedly that Brownian networks are potentially applicable as approximate models of more complex systems where, for example, processing activities involve simultaneous usage of several resources (servers) or require several different materials as inputs. Hereafter such physical systems will be referred to generically as "processing networks," or "stochastic processing networks."

Section 2 of this paper will describe in broad outline an extremely general family of stochastic processing networks, one which includes not only models with simultaneous resource requirements and multiple inputs to a single processing activity, but also alternative means of accomplishing a given task, and probabilistic work flow that may depend on which of those alternative means is chosen. Consider, for example, a manufacturing system where either a new machine or an old machine can be used to perform a certain operation, and suppose that there is a 50% chance that an additional "rework" operation will be required if the old machine is used, but no chance of that additional requirement if the new machine is used; further suppose that there is a choice as to which resources will be used to perform the rework operation if it is required. Readers will see that by properly defining "processing activities" and "job classes" (these are both primitive concepts in our general description of a processing network) it is relatively straightforward to represent such systems within the proposed modeling framework. The focus of this paper is on formulation and "soft analysis" of Brownian models for processing networks in that general family.

A recent paper [10] put forth a general explanation of the "state space collapse" that is a key to the tractability of Brownian networks. That is, the authors explained in general mathematical terms how a stochastic control problem associated with a Brownian network reduces to an "equivalent workload formulation" of lower dimension. The theory developed in [10] served to unify various ad hoc analyses of specially structured Brownian networks that had appeared over a span of years, and the authors also provided a general interpretation of the lower-dimensional workload formulation in terms of what they called "reversible displacements." That interpretation did not involve processing network language, invoking instead the physical image of a particle moving in high-dimensional space subject to both random Brownian displacements and purposeful control displacements.

This paper connects that general theory with the special structure one sees in Brownian models of processing networks, thereby extending the seminal work of Laws [14, 15] on Brownian models of alternate routing problems. To make that connection requires a modest amount of additional mathematics, and several questions left open in [10] will be answered along the way. To explain more precisely the character and contributions of this paper, it will be efficient to make a small diversion at this point, explaining what is meant by a Brownian network and summarizing the main result proved in [10].

Let $X=\{X(t), t\geq 0\}$ be an $m$-dimensional Brownian motion with respect to a given filtration on a fixed probability space. (A process that is adapted to the given filtration is nonanticipating with respect to $X$.) We denote by $\theta$ and $\Sigma$ the drift vector and covariance matrix respectively of $X$, assuming throughout the $X(0)=0$ almost surely. Also given are an $m\times n$ input–output matrix $R$, a $p\times n$ capacity consumption matrix $K$ and an initial inventory vector $q\in\mathbb{R}^m_+$. In [10] the drift vector of $X$ and the initial inventory vector were denoted by $\mu$ and $z$, respectively; with those exceptions, the notation in this paper agrees with that used in [10]. An admissible control is an $n$-dimensional process $Y=\{Y(t), t\geq 0\}$ such that

$$Y \text{ is adapted to the given filtration,} \tag{1.1}$$

$$U(\cdot) \text{ is nondecreasing with } U(0)\geq 0 \text{ and} \tag{1.2}$$

$$Z(t)\geq 0 \qquad \text{for all } t\geq 0, \text{ where} \tag{1.3}$$

$$Z(t)=q+X(t)+RY(t) \qquad \text{for all } t\geq 0 \text{ and} \tag{1.4}$$

$$U(t)=KY(t) \qquad \text{for all } t\geq 0. \tag{1.5}$$

The crucial element in the equivalent workload formulation developed in [10] is a matrix $M$ defined as follows. First, the space of "reversible displacements" referred to above is $\mathcal{N}=\{\delta\in\mathbb{R}^m : \delta=Ry, Ky=0\}$. Now let $\mathcal{M}$ be the orthogonal complement of $\mathcal{N}$, denote by $d$ the dimension of the linear space $\mathcal{M}$, and let $M$ be any $d\times m$ matrix whose rows are a basis for $\mathcal{M}$. That is, one can choose $M$ to be any $d\times m$ matrix with the following property:

$$\delta\in\mathcal{N} \quad \text{if and only if} \quad M\delta=0. \tag{1.6}$$

An equivalent characterization proved in [10] is the following: one can choose $M$ as any $d\times m$ matrix of full row dimension such that, for some $d\times p$ matrix $G$,

$$MR=GK. \tag{1.7}$$

The $m$-dimensional process $Z$ defined by (1.4) is generically called an inventory process, and it is $Z(t)$ that summarizes the state of the system at time $t$ in our original Brownian network model (1.1)–(1.5). The matrix $M$ is used to derive from $Z$ an associated "workload process" $W$ via

$$W(t) = MZ(t), \qquad t \geq 0, \tag{1.8}$$

and it is $W(t)$ that summarizes the state of the system at time $t$ in our equivalent workload formulation, or reduced Brownian network. That is, two state vectors $z$ and $z'$ are equivalent in the Brownian system model if their associated workload vectors $w = Mz$ and $w' = Mz'$ are identical. For purposes of this paper it is not even necessary to write out the equivalent workload formulation, but one point is crucial. The workload dimension $d$ always satisfies $d \leq m$, and $d$ is typically much smaller than $m$ in realistic applications, so adoption of the reduced state description (1.8) typically represents a substantial simplification of the original Brownian network. Moreover, as shown by Laws [14, 15] in his analysis of alternate routing problems, a careful examination of the "state space collapse" embodied in (1.8) may yield valuable qualitative information about essential system structure.

Once the Brownian network model (1.1)–(1.5) and its equivalent workload formulation have been described, certain questions naturally suggest themselves. To begin at the very beginning, *how are the data of a Brownian network, especially the matrices $R$ and $K$ from which $M$ is computed, derived from the data of a processing network whose behavior we wish to approximate*? That question was answered for certain classes of queueing networks, or at least a general recipe was proposed for those classes of queueing networks, in [6]. Extensions to more complex queueing networks that involve dynamic routing decisions have been suggested in a series of more recent papers surveyed by Kelly and Laws [13]. Here the question will be addressed in broad terms for the case of *open* processing networks, in which items or materials to be processed are generated exogenously, the network's resources or servers can undertake various activities to accomplish that processing and all materials eventually leave the system when their processing is complete. This paper's treatment of Brownian model formulations for open processing networks, which culminates in Section 5, is novel in two regards. First, it embraces the general model class referred to in the opening paragraph of this introduction, without reference to the special features that distinguish queueing network models. Second, there is explicit discussion of the need to distinguish between what are here called basic and nonbasic activities. In particular, when some of the available activities are nonbasic, in a precise sense explained in Section 2, the appropriate definition of $K$ for an approximating Brownian model involves augmentation of the original problem data. (In the notation used here, the original capacity consumption matrix for an open processing network is denoted by $A$, and rows are added to it to form $K$. The original matrix $A$ is nonnegative, but the rows that are added to it contain negative elements.) The existing literature on Brownian network approximations contains no mention of nonbasic activities, although the subject is potentially important and somewhat subtle.

A second basic modeling question, inextricably intertwined with the first, is *how to articulate the "heavy traffic" condition* required to justify a Brownian approximation. We address that matter in Section 2 by means of linear programming, extending in an obvious way the analysis of alternate routing problems by Laws [14, 15]. The linear program (LP) used to define heavy traffic uses only first-order data (average arrival rates, average service rates and so forth), and it is called *the static planning problem* in this paper. Its dual, whose basic feasible solutions define what Laws calls "generalized cut constraints," is here called the *workload definition problem.*

Perhaps the most obvious question left unanswered in [10] is the following: *Is there a canonical choice of the basis matrix $M$*, used to define the workload process $W$ in (1.8), *that is natural in a processing network context*? For Brownian models of open processing networks, a qualified affirmative answer is provided by the development in Sections 2–5, which extends the analysis by Laws [14, 15] of a more restrictive model class. We define a "canonical choice" of the matrix $M$ such that each row corresponds to a generalized cut constraint which is binding at the optimal solution of the static planning problem. (This choice of $M$ is not necessarily unique, so it produces *a* canonical representation of workload rather than *the* canonical representation.) As a by-product of this argument we obtain an interpretation of $\mathcal{M}$ in terms of a fluid model of the original processing network. As explained in Section 5, the deterministic fluid model provides a useful complement to the Brownian model that is our central focus.

In the context of open processing networks, two other natural questions about the matrix $M$ concern its nonnegativity and factorability. With regard to the former, a monotonicity condition on the original problem data is presented in Section 6, and it is shown to imply that our canonical choice of $M$ is nonnegative. In general, however, there need not exist a nonnegative basis matrix for the space $\mathcal{M}$. Another natural condition is presented in Section 6 which ensures that our canonical choice of $M$ admits a factorization related to the notion of resource pooling. Broadly speaking, resource pooling occurs when a system manager has alternative means available for processing a given set of inputs, and as Kelly and Laws [13] have rightly emphasized, factorizations of the sort discussed in Section 6 yield insights that are important in every application domain.

Careful study of examples and special model structures is crucial for understanding the general phenomenon of state space collapse in Brownian networks. Several examples are discussed in Sections 6 and 7 to illustrate specific points, but to develop a full appreciation for the specially structured model classes that motivate our general theory, readers who are new to the subject should review earlier work. The historical progression begins in [6] with analysis of open multiclass networks where the only mode of dynamic control involves sequencing decisions. The original dimension of the Brownian model for such a network equals the number of customer classes, but the equivalent workload formulation has dimension equal to the number of servers. Brownian models of closed multiclass queueing networks are also discussed in [6] and in [4], but the most comprehensive treatment of those models and their equivalent workload formulations is contained in Section 6 of [10]. The most surprising and important results regarding state space collapse arise in the context of multiclass queueing networks with alternate routing decisions, in addition to dynamic sequencing decisions: the Brownian models of such networks typically have equivalent workload formulations with dimension strictly less than the number of servers. The survey paper by Kelly and Laws [13] describes a series of illuminating examples, one of which will be reviewed in Section 7 of this paper. Section 3 of [10] analyzes another example of similar character.

All of the examples considered in Section 6 of this paper involve processing activities with multiple inputs, so they are not queueing network models as that phrase is normally interpreted, and most have been chosen to illustrate some phenomenon that cannot occur in a conventional queueing network context. For one of our examples there exists no choice of the basis matrix $M$ that has all elements nonnegative. Another example has an equivalent workload formulation whose dimension $d$ is strictly larger than the number of servers, and that same example illustrates the role of nonbasic activities in a Brownian network model. The examples described in Section 6 also show how a given set of processing capabilities may yield very different Brownian network models depending on the nature of the input streams to be processed.

The queueing network model considered in Section 7 of this paper, which originated in the Ph.D. dissertation of Laws [14], provides an example where our "canonical choice" of the basis matrix $M$ is nonunique. (In general, however, the approach developed here does narrow the canonical choice to finitely many candidates.) Section 7 not only describes that specific example but also explains how the theory developed in this paper relates to the earlier work of Kelly and Laws [13, 14, 15].

In concluding this introduction, it is appropriate to emphasize two problems or issues *not* addressed in the paper. First, open processing networks are described only in broad terms, without any attempt at precise mathematical formulation of the general model class. To develop such a formulation is a major undertaking in itself, and here we provide only enough detail as to make plausible the proposed Brownian approximation, fully acknowledging that rigorous justification of the Brownian model as a heavy traffic approximation may depend on modeling choices that lie below the level of resolution in our treatment. The other problem or issue referred to above is how to translate an optimal control policy for the approximating Brownian model, assuming one can be found, into a nearly optimal policy for the original processing network. That general problem, which is known to be difficult and subtle, is a focus of current research [8, 9] but simply will not be broached in this paper.

**2. The static planning problem for an open network.** Consider an open processing network with $r$ servers (or processing resources), $m$ materials (or stocks, or job classes) and $n$ processing activities. It will be convenient to imagine that each material is stored in a dedicated buffer, and the terms "buffer" and "inventory" will occasionally be used as synonyms for "material." Each activity requires certain materials as inputs, and it may either destroy those inputs or produce other materials as outputs. Let us denote by $R_{ij}$ the average amount of material $i$ consumed per unit of activity $j$, with a negative value interpreted to mean that activity $j$ is a net producer of material $i$. It is also the case that activities consume the capacities of associated resources. Assuming as a matter of convention that each resource has available one unit of capacity per unit of time, let $A_{kj}$ be the average amount of resource $k$ capacity consumed per unit of activity $j$. In general, material consumption, material production and capacity consumption could all be stochastic, depending on how units of activity are defined, but we shall deal only with average rates for the time being. Similarly, let $\lambda_i$ be the average rate at which material $i$ is exogenously generated, assuming $\lambda_i > 0$ for at least one $i$. The $m \times n$ input–output matrix $R$, the $r \times n$ nonnegative capacity consumption matrix $A$, and the nonnegative $m$-dimensional column vector $\lambda$ provide the first-order data for our open processing network.

The static planning problem referred to in Section 1 is the following: choose a scalar $\rho$ and an $n$-dimensional column vector $x$ so as to

$$\text{minimize} \quad \rho \tag{2.1}$$

$$\text{subject to} \quad Rx = \lambda, \qquad Ax \leq e\rho \quad \text{and} \quad x \geq 0, \tag{2.2}$$

where $e$ is the $r$-vector of ones. One interprets $x_j$ as the average rate at which activity $j$ is undertaken, expressed in units of activity per unit of time, and $\rho$ as an upper bound on the utilization rates for our various resources under the processing plan $x$. In the static planning problem (2.1) and (2.2), one seeks to minimize the maximal utilization rate $p$ subject to the requirement that average activity rates be nonnegative and that exogenously generated inputs be processed to completion without other inventories being generated ($Rx = \lambda$).

A natural first question to ask is whether there exists an $(x,\rho)$ satisfying (2.2) with $\rho \le 1$, and in that regard it is useful to consider the dual linear program: choose an $m$-dimensional row vector $\mu$ and an $r$-dimensional row vector $\pi$ so as to

$$\text{maximize} \quad \mu\lambda \tag{2.3}$$

$$\text{subject to} \quad \mu R \le \pi A, \qquad \pi e = 1 \quad \text{and} \quad \pi \ge 0. \tag{2.4}$$

Actually, the constraint $\pi e = 1$ in (2.4) is first expressed with a $\le$ inequality, but that can obviously be converted to an equality as shown, because the capacity consumption matrix $A$ is nonnegative. One can interpret $\mu_i$ as the "total work content" attributed to a unit of material $i$, and $\pi_k$ as the relative capacity of server $k$. The constraint $\pi e = 1$ sets total system capacity to one as a matter of convention, which justifies our description of $\pi_k$ as a "relative" capacity. The constraint $\mu R \le \pi A$ in (2.4) demands that workload contributions $\mu_i$ be attributed to the various materials $i = 1, \ldots, m$ and relative processing capacities $\pi_k$ be attributed to the various resources $k = 1, \ldots, r$ in such a way that no activity accomplishes more "work" than the total capacity it consumes. Subject to that crucial constraint, one strives to maximize $\mu\lambda$, which is the total amount of work attributed to the materials that are exogenously generated per unit of time. These interpretations, which are discussed further in Section 4, justify the description of (2.3) and (2.4) as a "workload definition problem."

Let us denote by $\mathcal{D}$ (mnemonic for *dual)* the set of all $(\mu,\pi)$ pairs satisfying (2.4), and let $\mathcal{E} = \{(\mu^l, \pi^l) : l = 1, \ldots, L\}$ be the extreme points of the polytope $\mathcal{D}$. In other words, $\mathcal{E}$ consists of the $L$ distinct *basic feasible solutions* for our workload definition problem (2.3) and (2.4). The fundamental theorem of linear programming tells us that the minimal objective value $\rho^*$ for the primal problem (2.1) and (2.2) equals the maximal value of $\mu\lambda$ over all $(\mu,\pi) \in \mathcal{D}$, or equivalently, equals the maximum of $\mu^l\lambda$ over all extreme points $l = 1, \ldots, L$. Thus we have the following.

PROPOSITION 1. *$\rho^* \le 1$ if and only if*

$$\mu^l \lambda \le 1 \qquad \textit{for all } l = 1, \ldots, L. \tag{2.5}$$

Inequality (2.5) corresponds to what Kelly and Laws [13, 14, 15] called "generalized cut constraints," that language reflecting their interest in multicommodity network flows. The "heavy traffic condition" for our open processing network is now articulated as follows.

ASSUMPTION 1. The static planning problem (2.1) and (2.2) has a unique optimal solution $(\rho^*, x^*)$. Moreover, that solution has $\rho^* = 1$ and $Ax^* = e$.

This assumption will remain in force throughout the paper, and several of its aspects require comment. First, given the vector $\lambda$ of external arrival rates, Assumption 1 says that every server must be fully utilized ($Ax^* = e$) if we are to avoid inventory build-ups over time. This situation where *all* servers are "critically loaded" may seem rather special, but the assumption really means that servers which need not be fully utilized to handle the given load have simply been deleted from our model. That is, the model under discussion provides explicit representation of just the "bottleneck subnetwork" composed of $r$ critically loaded servers. (Of course, one would like to see a formal proof that the influence of noncritical resources is negligible in some sense, but that matter will not be dealt with here.) Also, it might seem unlikely that more than one resource would be critically loaded, but multiple bottlenecks are commonplace in manufacturing systems because of capacity balancing, and Laws [15] has shown how alternate routing capabilities lead to "resource pools" whose constituent servers all approach the critical loading regime together. Next, readers might object that "critical loading" is equated with the case $\rho^* = 1$ in this discussion, rather than the more generous interpretation to mean $\rho^*$ *near* 1. In that regard Assumption 1 should not be taken quite literally: as is standard in the theory of diffusion approximations, one can consider a sequence of systems whose first-order data ($R$, $A$ and $\lambda$) approach limits, requiring only that the *limiting data* satisfy Assumption 1.

Finally, it may seem needlessly restrictive to require that the static planning problem have a unique solution, but multiple optima (with all servers fully utilized in every optimal solution) lead to Brownian approximations of a more complex form than considered in this paper, and all examples that have been analyzed in the literature to date satisfy our assumption. For open queueing networks with alternate routing, the uniqueness assumption amounts to the following: there is just one way of splitting arrivals in each external input stream among available alternate routes so that the no server is loaded beyond its capacity by the resultant flows.

Hereafter, the vector $x^*$ of average activity levels will be said to constitute a *nominal processing plan.* Because of stochastic variability in exogenous inputs and endogenous processing, actual average activity rates over moderate time spans may differ from these nominal rates, but the differences must be small over long time spans if excessive inventories are to be avoided.

Given that $\rho^* = 1$ by Assumption 1, it follows from Proposition 1 and the remarks immediately preceding it that $\mu^l \lambda = 1$ for at least one $l$. That is, at least one of the generalized cut constraints (2.5) must be binding, and hereafter we denote by $L^*$ the number of such binding constraints. Next let $d$ be the dimension of the linear space spanned by $\{\mu^l : 1 \le l \le L^*\}$. (Eventually, it will be shown that this number $d$ equals the "effective system dimension" for an appropriate Brownian approximation of the open processing network, so there is no conflict with the notation of Section 1.) Without loss of generality, we assume that the extreme points of $\mathcal{D}$ are numbered so that

$$\mu^l \lambda = 1 \qquad \text{for } 1 \le l \le L^*, \qquad \text{but } \mu^l \lambda < 1 \qquad \text{for } L^* < l \le L, \tag{2.6}$$

and so that $\mu^1, \ldots, \mu^d$ are also linearly independent. Finally, for future purposes we define a $d \times m$ matrix $M$ and a $d \times r$ matrix $\Pi$ via

$$M = \begin{bmatrix} \mu^1 \\ \vdots \\ \mu^d \end{bmatrix} \quad \text{and} \quad \Pi = \begin{bmatrix} \pi^1 \\ \vdots \\ \pi^d \end{bmatrix}. \tag{2.7}$$

This matrix $M$ is the canonical basis matrix referred to in Section 1; its use in that role will be justified in the following sections. (Obviously, our "canonical choice" of $M$ is unique if and only if $d = L^*$.)

Let us denote by $b$ the number of activities $j$ such that $x_j^* > 0$, calling these *basic* activities, and let activities be numbered so that the basic ones are $1, \ldots, b$. Activities $b+1, \ldots, n$ will be called *nonbasic*, and it will be convenient for later purposes to partition $R$ and $A$ in the following way:

$$R = [H \quad J] \quad \text{and} \quad A = [B \quad N], \tag{2.8}$$

where $H$ and $B$ both have $b$ columns. Thus $H$ and $B$ are the submatrices of $R$ and $A$, respectively, that correspond to basic activities. Because the optimal solution $(\rho^*, x^*)$ of our static planning problem (2.1) and (2.2) may be degenerate, an optimal basis for that linear program may include one or more of the variables $x_{b+1}, \ldots, x_n$ but that need not concern us here. What *is* important for our purposes is that each of the pairs $(\mu^1, \pi^1), \ldots, (\mu^d, \pi^d)$ is an optimal solution for the dual linear program (2.3) and (2.4), and thus $(\mu^l H = \pi^l B$ for $l = 1, \ldots, d$ by complementary slackness. (I.e., if any component of $\mu^l H$ were strictly smaller than the corresponding component of $\pi^l B$, one would have $\mu^l \lambda < \rho^* = 1$, which is a contradiction.) Given the definitions (2.7), one can write this in matrix form as

$$MH = \Pi B. \tag{2.9}$$

**3. A balanced fluid model of dynamic flow management.** In this section we consider a deterministic fluid model of the open processing network described in Section 2, again using only the first-order data $R$, $A$ and $\lambda$. One may paraphrase Assumption 1 by saying that there is a perfect balance between the vector $\lambda$ of external input rates and the processing capacities of the system's $r$ servers, so we are examining a *balanced fluid model* of the sort emphasized in [7], whereas fluid model analyses of stochastic control systems that have been published to date [1, 2, 3, 5, 20] have been concerned primarily with the case where available capacity strictly exceeds what is needed to process external inputs. Readers will eventually see that the two formal propositions proved in this section could easily be framed without any mention of a fluid model, but developing them in this way improves one's intuitive understanding of state space collapse in Brownian networks (see Section 5).

In the fluid model a policy for dynamic flow management takes the form of a measurable function $\alpha : [0,\infty) \to \mathbb{R}^n_+$, where $\alpha_j(t)$ is interpreted as the rate at which activity $j$ is undertaken at time $t$(expressed in units of activity per unit of time). Assuming initial inventory vector $q \in \mathbb{R}^m_+$, the inventory trajectory $z(\cdot)$ generated from $\alpha(\cdot)$ is

$$z(t) = q + \lambda t - R \int_0^t \alpha(s)\, ds, \qquad t \geq 0. \tag{3.1}$$

The policy $\alpha(\cdot)$ is called *feasible* if it satisfies $\alpha(t) \geq 0$ for all $t \geq 0$, $z(t) \geq 0$ for all $t \geq 0$, and

$$A\alpha(t) \leq e \qquad \text{for all } t \geq 0. \tag{3.2}$$

A state $q' \in \mathbb{R}^m_+$ is said to be *reachable* from $q$ if there exists a feasible policy $\alpha(\cdot)$ and a time $t > 0$ such that $z(t) = q'$, where $z(\cdot)$ is defined by (3.1). We say that states $q$ and $q'$ *communicate* if each is reachable from the other.

PROPOSITION 2. *Suppose $q$, $q' \in \mathbb{R}^m_+$ and let $\delta = q' - q$. The following are equivalent:*

(i) *$q'$ is reachable from $q$.*
(ii) *There exists $x \geq 0$ such that $Ax \leq e$ and $Rx = \lambda - (1/t)\delta$ for some $t > 0$.*
(iii) *$\mu^l \delta \geq 0$ for all $l = 1, \dots, L^*$.*

PROOF. It follows immediately from Proposition 1 that (ii) holds if and only if, for each $l = 1, \dots, L$, one has $\mu^l(\lambda - (1/t)\delta) \leq 1$ for some $t > 0$. By (2.6) the latter condition holds if and only if $\mu^l \delta \geq 0$ for all $l = 1, \dots, L^*$, which is (iii). Thus we have proved that (ii) and (iii) are equivalent. Now suppose that (i) holds, which means that

$$q' = q + \lambda t - R \int_0^t \alpha(s)\, ds \tag{3.3}$$

for some feasible policy $\alpha(\cdot)$ and some $t > 0$. Defining $x = \int_0^t \alpha(s)\, ds/t$, the feasibility of $\alpha(\cdot)$ implies $x \geq 0$ and $Ax \leq e$, so (ii) is satisfied by this particular $x$ and $t$. Conversely, if (ii) holds we can take $\alpha(s) = x$ for $0 \leq s \leq t$ to satisfy (i). Thus (i) and (ii) are equivalent, which completes the proof. $\square$

Anticipating the treatment of Brownian models to follow in Section 5, it will be useful to re-express fluid model dynamics in terms of a centered cumulative time allocation $\beta$, as follows. Given any feasible policy $\alpha(\cdot)$, one defines an associated control $\beta$ via

$$\beta(t) = x^* t - \int_0^t \alpha(s)\, ds, \qquad t \geq 0. \tag{3.4}$$

Thus components of $\beta(t)$ represent cumulative activity levels over the time interval $[0, t]$, expressed as decrements from the vector $x^* t$ of nominal activity levels for that interval. Recalling that $Rx^* = \lambda$, we can rewrite the basic system equation (3.1) as

$$z(t) = q + R\beta(t), \qquad t \geq 0. \tag{3.5}$$

A feasible control $\beta$ must have $z(\cdot) \geq 0$, of course, and to express other policy constraints in terms of $\beta$ it will be convenient to set

$$p = r + n - b \tag{3.6}$$

and then define a $p \times n$ matrix $K$ as follows:

$$K = \begin{bmatrix} B & N \\ 0 & -I \end{bmatrix}. \tag{3.7}$$

Comparing (3.7) and (2.8), one sees that the first $r$ rows of $K$ are the capacity consumption matrix $A$, and the negative identity matrix appearing in (3.7) is of dimension $n - b$, which is the number of nonbasic activities in our static planning problem.

Returning now to the reformulation of our fluid model in terms of $\beta$, let us define

$$u(t) = K\beta(t), \qquad t \geq 0 \tag{3.8}$$

and consider the following policy constraint:

$$u(\cdot) \text{ is nondecreasing with } u(0) = 0. \tag{3.9}$$

Recalling that $Ax^* = e$ by Assumption 1 (i.e., all servers are loaded to full capacity under the nominal processing plan $x^*$), one sees from (3.4), (3.7) and (3.8) that the first $r$ components of constraint (3.9) are equivalent to (3.2), while the last $p - r$ components reexpress the requirement that $\alpha_j(\cdot) \geq 0$ for $j = b+1, \ldots, n$ (these $j$ correspond to nonbasic activities). The requirement that $\alpha_j(\cdot) \geq 0$ for $j = 1, \ldots, b$ can be reexpressed as follows:

$$\beta_j(t) - \beta_j(s) \leq (t - s)x_j^* \qquad \text{for } j = 1, \ldots, b \text{ and } 0 \leq s \leq t. \tag{3.10}$$

To summarize, a control $\beta : [0, \infty) \to \mathbb{R}^n$ is deemed feasible in our reformulated fluid model if it is measurable and the following hold: the process $z$ defined by (3.5) satisfies $z(t) \geq 0$ for all $t \geq 0$; the process $u$ defined by (3.8) satisfies (3.9) and $\beta$ further satisfies (3.10). Having introduced the matrix $K$ in conjunction with the reformulation, we can now state the main mathematical result of this paper.

PROPOSITION 3. *Suppose that* $q$, $q' \in \mathbb{R}^m_+$ *and let* $\delta = q' - q$. *The following are equivalent:*

(i) $q$ *and* $q'$ *communicate.*
(ii) *There exists* $y \in \mathbb{R}^n$ *such that* $Ry = \delta$ *and* $Ky = 0$.
(iii) $M\delta = 0$ *(i.e.,* $Mq' = Mq$*).*

PROOF. It follows from Proposition 2 that $q$ and $q'$ communicate if and only if $\mu^l \delta = 0$ for all $l = 1, \dots, L^*$. The rows of $M$ have been chosen so that the latter condition holds if and only if $M\delta = 0$. Thus (i) and (iii) are equivalent.

Now suppose that (ii) holds, and define $x = x^* - (1/t)y$, where $t$ is large enough to ensure $x_j \geq 0$ for all $j = 1, \dots, b$ (recall that $x_j^* > 0$ for all such $j$). Because $Ky = 0$, we have from the definition (3.7) of $K$ that $Ay = 0$ and $y_j = 0$ for $j = b+1, \dots, n$. Recalling that $Ax^* = e$, $Rx^* = \lambda$ and $x_j^* = 0$ for $j = b+1, \dots, n$, one then has that $Ax = e$, $Rx = \lambda - (1/t)\delta$ and $x \geq 0$. Thus $q'$ is reachable from $q$ by Proposition 2. Now one can use that same argument with $-y$ in place of $y$ to prove that $q$ is reachable from $q'$, so (ii) implies (i).

Finally, suppose that (i) holds. Without loss of generality we can assume the existence of a $t > 0$ such that $q'$ is reachable from $q$ in exactly $t$ time units and $q$ is also reachable from $q'$ in exactly $t$ time units. [Having reached a target state $q^*$ at some time $\tau$, one can always maintain $z(t) = q^*$ for $t \geq \tau$ by simply taking $\alpha(t) = x^*$ for $t \geq \tau$.] Arguing exactly as in the proof of Proposition 2, we then have the following: there exist vectors $x$, $x' \geq 0$ such that $Ax \leq e$, $Ax' \leq e$, $Rx = \lambda - (1/t)\delta$ and $Rx' = \lambda + (l/t)\delta$. Defining $x'' = \frac{1}{2}(x + x')$, one then has $Ax'' \leq e$, $Rx'' = \lambda$ and $x'' \geq 0$. Thus one obtains a feasible solution for the static planning problem (2.1) and (2.2) by taking $x = x''$ and $\rho = 1$. But Assumption 1 says that the unique optimal solution $(\rho^*, x^*)$ has $\rho^* = 1$, so it must be that $x'' = x^*$. Then it must be that $Ax = Ax' = e$ and $x_j = x'_j = 0$ for $j = b+1, \dots, n$ because $Ax^* = e$ and $x_j^* = 0$ for $j = b+1, \dots, n$. Thus we can satisfy (ii) by taking $y = t(x^* - x)$. This shows that (i) implies (ii), completing the proof. □

As the notation suggests, the matrix $K$ defined by (3.7) will serve as the capacity consumption matrix in our Brownian network model (see Section 4). The last order of business in this section is to exhibit a $d \times p$ matrix $G$ such that $R$, $G$, $K$ and the canonical basis matrix $M$ in (2.7) jointly satisfy the definitive relationship (1.7). The natural choice is

$$G = [\Pi \quad \Lambda] \qquad \text{where } \Lambda = \Pi N - MJ \geq 0. \tag{3.11}$$

The nonnegativity of $\Lambda$ declared in (3.11) is established as follows: first, each row of $\Lambda$ corresponds to a pair $(\mu, \pi)$ that is a feasible solution of the dual linear program (2.3) and (2.4); a constraint of that LP is $\mu R \leq \pi A$, which implies $\mu J \leq \pi N$ by (2.8); thus $\pi N - \mu J \geq 0$, implying $\Pi N \geq MJ$ as claimed. Combining (3.11) with (2.8), (2.9) and the definition (3.7) of $K$, one has that

$$G \geq 0 \quad \text{and} \quad MR = GK. \tag{3.12}$$

## 4. Further interpretation of the dual variables $\boldsymbol{\mu}$ and $\boldsymbol{\pi}$.

In Section 2 it was suggested that the variable $\mu_i$ appearing in our dual linear program (2.3) and (2.4) be interpreted as the workload contribution, or total work content, per class $i$ job, and that $\pi_k$ be interpreted as the relative capacity of server $k$. To sharpen those interpretations, we now restate the primal linear program (2.1) and (2.2) as follows: given an arbitrary $m$-vector $q$, find a scalar $\tau$ and $n$-vector $x$ so as to

$$\text{minimize} \quad \tau \tag{4.1}$$

$$\text{subject to} \quad Rx = q, \qquad Ax \leq e\tau \quad \text{and} \quad x \geq 0. \tag{4.2}$$

Adopting the deterministic fluid model of system dynamics that was described in Section 3, we interpret $\tau$ as the length of a time interval during which exogenous inputs are to be turned off, while $q_i$ is interpreted as the amount of material which must be removed from buffer $i$ during that interval (if $q_i$ is negative, this means that processing activities during the interval must *add* material to buffer $i$). The decision variable $x_j$ is interpreted as the total amount of activity $j$ to be undertaken during the interval $[0, \tau]$, and by expressing capacity constraints in the form $Ax \leq e\tau$, we are effectively restricting attention to controls $\alpha$ (see Section 3) in which each activity $j$ is undertaken at constant rate $\alpha_j(t) = x_j/\tau$ $(0 \leq t \leq \tau)$. Given that our objective (4.1) is to minimize the length $\tau$ of the processing interval, it is easy to verify that the restriction to constant processing rates actually involves no loss of generality.

The dual of the least-time control problem (4.1) and (4.2) is exactly as stated earlier in (2.3) and (2.4), except that $q$ is substituted for $\lambda$ in the objective (2.3). We call the optimal objective value $\tau^*$ the *minimum time to execute* the target vector $q$ of buffer content changes, abbreviating this as MTTE. If all components of $q$ are nonnegative, as occurs when $q = \lambda$, then one can think of MTTE as an acronym for *minimum time to emptiness*, starting with the buffer contents vector $q$. When it is desirable to emphasize the dependence of $\tau^*$ on $q$ we write $\tau^*(q)$.

Returning to the case $q = \lambda$ that was considered in Section 2, recall that Assumption 1 includes the condition $\rho^* = \tau^*(\lambda) = 1$. From this it follows easily that $\tau^*(t\lambda) = t$ for any $t > 0$. For purposes of interpretation, let us assume initially that there exists a *unique* optimal dual solution $(\mu^*, \pi^*)$ when $q = \lambda$. Then that same pair $(\mu^*, \pi^*)$ is uniquely optimal for the dual problem when $q = t\lambda$, where $t > 0$ is arbitrary. Uniqueness of the optimal dual solution implies that $\tau^*(t\lambda+\delta)=\tau^*(t\lambda)+\mu^*\delta=t+\mu^*\delta$ for any $m$-vector $\delta$ and all $t$ sufficiently large. That is, the minimum time to emptiness when starting with exactly $t$ time units of exogenous input flow in each buffer is itself $t$, and $\mu_i^*$ is the rate at which MTTE increases as the initial content of buffer $i$ increases. Thus, equating "system workload" with minimum time to emptiness, one can accurately say that $\mu_i^*$ represents the workload contribution per incremental unit of material in buffer $i$, starting from a base case where the vector of initial buffer contents is $q=t\lambda$ with $t$ large.

Continuing to assume that the optimal solution $(\mu^*,\pi^*)$ of (2.3) and (2.4) is unique, one can interpret the optimal dual variable $\pi_k^*$ as the rate at which MTTE would decrease from the initial value $\tau^*(t\lambda)=t$ if additional capacity were available from server $k$. That is, if server $k$ were available for $\varepsilon$ time units of "preprocessing" before the clock measuring time-to-emptiness began to run, and if the system manager used that time optimally, then MTTE would decrease by $\pi_k^*\varepsilon$. Conversely, if server $k$ were forced to remain idle during the first $\varepsilon$ time units while other servers were working, then MTTE would increase by $\pi_k^*\varepsilon$. With this interpretation, the constraint $\pi e=1$ in (2.4) is entirely natural, because making all servers available for $\varepsilon$ time units of preprocessing would obviously decrease MTTE by $\varepsilon$. Using the language of fluid models and again equating "work in the system" with MTTE, one can describe $\pi_k^*$ as the rate at which server $k$ alone "drains work from the system," starting from a vector of initial buffer contents (or more generally, a target vector $q$ of buffer content changes) which is sufficiently close to $\lambda t$ for some $t>0$.

Suppose now that the dual linear program, (2.3) and (2.4) has two or more basic optimal solutions $\{(\mu^l,\pi^l): 1\le l\le L^*\}$. Then one has

$$\tau(t\lambda+\delta)=t+\max_{1\le 1\le L^*}\mu^l\delta \tag{4.3}$$

for any $m$-vector $\delta$ and $t$ sufficiently large. Consider a target vector $q$ of buffer content changes having the form $q=t\lambda+\delta$, where $t$ is large and $\delta$ is such that the maximum in (4.3) is uniquely achieved by one particular $l$. Then the basis associated with $(u^l,\pi^l)$ will be uniquely optimal for our least-time control problem (4.1) and (4.2), and the optimal dual variables $u_i^l$ and $\pi_k^l$ will have exactly the interpretations given above. It is noteworthy that some servers $k$ may have $\pi_k^l=0$, which means that such servers are rendered *noncritical* by a perturbation in direction $\delta$: that is, MTTE is not decreased by additional server $k$ capacity if the target vector of buffer content changes is $q=t\lambda+\delta$. In the same way, one may have $u_i^l=0$ for some materials $i$, which means that incremental units of material $i$ have no effect on MTTE after a perturbation in direction $\delta$. (This may occur, e.g., if such incremental material can be processed using only servers which have been rendered noncritical by the perturbation.) For an illustration of these phenomena, see Section 7.

**5. A Brownian model of dynamic flow management.** Recall that in Section 1, given a Brownian network with data $R$ and $K$, we defined the space of reversible displacements $\mathcal{N} = \{\delta \in \mathbb{R}^m : \delta = Ry, Ky = 0\}$. A basis matrix for the orthogonal space $\mathcal{M}$ is any $M$ satisfying (1.6). Comparing (1.6) with Proposition 3, and noting that any $\delta \in \mathbb{R}^m$ can be written as the difference of two nonnegative vectors $q$ and $q'$, we have the following.

PROPOSITION 4. *Suppose that the matrices $R$ and $K$ of a Brownian network are derived from the first-order data ($R$, $A$ and $\lambda$) of an open processing network satisfying Assumption 1, as described in Sections 2 and 3. Then the matrix $M$ defined in Section 2 is indeed a basis matrix for the space $\mathcal{M}$, and hence one feasible representation of the workload process $W(t)$ discussed in Section 1 is $W(t) = MZ(t)$.*

The remainder of this section is devoted to the following questions. First, why is it appropriate to define the data $R$ and $K$ of an approximating Brownian network in the way described here? Second, how should one define the drift vector $\theta$ and covariance matrix $\Sigma$ for the Brownian approximation in terms of processing network data? Finally, under what rescaling of time and state space would one expect to obtain convergence to the proposed Brownian approximation, and what rescaling would plausibly yield convergence to the deterministic fluid model described in Section 3?

To address these questions one must provide some detail about the stochastic microstructure of the open processing network, and one potential approach is the following. Generalizing the treatment of multiclass queueing networks in [6], one may take as primitive a collection of independent $m$-dimensional stochastic processes $F^j = \{F^j(t), t \geq 0\}$ indexed by $j = 0, 1, \ldots, n$. The $i$th component of the vector process $F^j(t)$ is denoted $F_i^j(t)$. We interpret $F_i^0(t)$ as the cumulative exogenous input of material $i$ up to time $t$. For $i = 1, \ldots, m$ and $j = 1, \ldots, n$ we interpret $F_i^j(t)$ as the cumulative amount of material $i$ consumed in the first $t$ units of activity $j$ undertaken, with a negative value indicating net production of the material rather than net consumption. A dynamic flow management policy takes the form of an $n$-dimensional stochastic process $T = \{T(t), t \geq 0\}$ with components $T_j(t)$. We interpret $T_j(t)$ as the cumulative amount of activity $j$ undertaken up to time $t$, and so the $m$-dimensional material inventory process, or buffer contents process, is given by

$$Q(t) = Q(0) + F^0(t) - \sum_{j=1}^{n} F^j(T_j(t)), \qquad t \geq 0. \tag{5.1}$$

Of course, the chosen process $T$ must satisfy a variety of constraints, and those will be discussed shortly. Making connection with the first-order data used in previous sections, we assume that $E[F^0(t)] \sim \lambda t$ as $t \to \infty$, and that $E[F^j(t)] \sim R^j t$ as $t \to \infty$ for $j = 1, \ldots, n$, where $R^j$ is the $j$th column of $R$. Setting $R^0 = \lambda$ to unify notation, the following sort of functional central limit theorem (FCLT) must also be assumed to justify the Brownian system model: for each $j = 0, 1, \ldots, n$ there exists an $m \times m$ covariance matrix $\Gamma^j$ such that the normalized vector process

$$\varepsilon[F^j(t/\varepsilon^2) - R^j t/\varepsilon^2], \qquad t \geq 0, \tag{5.2}$$

converges weakly as $\varepsilon \searrow 0$ to an $m$-dimensional Brownian motion with zero drift and covariance matrix $\Gamma^j$.

Suppose that the processing network satisfies Assumption 1, and define the *nominal inventory process*

$$\xi(t) = F^0(t) - \sum_{j=1}^{n} F^j(x_j^* t), \qquad t \geq 0. \tag{5.3}$$

From the FCLTs assumed for the individual flow processes $F^j$ one has the following (recall that $Rx^* = \lambda$): for small values of $\varepsilon > 0$, the scaled nominal inventory process

$$\varepsilon\xi(t/\varepsilon^2), \qquad t \geq 0, \tag{5.4}$$

is well approximated by a Brownian motion $X = \{X(t), t \geq 0\}$ with zero drift and covariance matrix

$$\Sigma = \Gamma^0 + x_1^*\Gamma^1 + \cdots + x_n^*\Gamma^n. \tag{5.5}$$

One obvious constraint on the chosen policy $T$ is that

$$Q(t) \geq 0, \qquad t \geq 0, \tag{5.6}$$

which we interpret to mean that an activity cannot be undertaken unless all of the required inputs are available. A proper articulation of the material availability constraint may actually require further restrictions on $T$, beyond (5.6), but such added restrictions are typically insignificant under the diffusion scaling to be imposed shortly. The simplest possible form that our capacity constraints can take is the following:

$$A[T(t) - T(s)] \leq (t - s)e \qquad \text{for } 0 \leq s \leq t. \tag{5.7}$$

Capacity constraints of this form are appropriate if capacity consumption is deterministic and occurs at uniform rates for all activities, and if processing resources can be shared by various activities in arbitrary proportions. Typically one expects the situation to be more complicated. For example, one resource might be a human operator who sets up different machines in a manufacturing facility so that they can run unattended after the operator has moved on to other tasks. Defining units of activity as total machine hours devoted to a given type of production, including both set-up time and run time, it is typically quite easy to determine the average number of operator hours required per unit of activity. But the precise articulation of operator availability and machine availability constraints is more complicated than (5.7) because all of an operator's time for a given activity is spent during just one part of a repeating production cycle, and an operator cannot split his or her time between set-up tasks at different machines. That is, a precise articulation of resource availability constraints may involve microcoordination concerns which are ignored in (5.7), but for current purposes we shall simply assume that (5.7) is an acceptable idealization under the rescaling that leads to a Brownian approximation.

Following the pattern established in Section 3 of [6], let us consider an arbitrary policy $T$ and define an associated $n$-dimensional centered process $V = \{V(t), t \geq 0\}$ via

$$V(t) = x^* t - T(t), \qquad t \geq 0. \tag{5.8}$$

The $j$th component $V_j(t)$ expresses the cumulative activity level for activity $j$ as a decrement from the nominal activity level $x_j^* t$, and under our heavy traffic assumption, only small values of this deviation control (relative to $t$) are of interest over long time spans. Using the identity $Ax^* = e$, one can reexpress the capacity constraints (5.7) as follows: defining $I(t) = AV(t)$ for all $t \geq 0$, one must choose the deviation control $V$ so that $I(\cdot)$ is nondecreasing with $I(0) = 0$. One interprets $I_k(t)$ as the cumulative amount of capacity for server $k$ that goes unused up to time $t$ (in a queueing network context the letter $I$ is mnemonic for *idleness*). As in the fluid model development in Section 3, it will compactify our mathematical theory to extend the $r$-dimensional "unused capacity process" $I$ as follows. Recalling the definition (3.7) of $K$, let

$$I(t) = KV(t), \qquad t \geq 0, \tag{5.9}$$

and then require that $V$ be chosen so that

$$I(\cdot) \text{ is nondecreasing with } I(0) = 0. \tag{5.10}$$

The first $r$ components of the process $I(\cdot)$ defined in (5.9) are the unused capacity process $AV(\cdot)$ discussed immediately above, so the first $r$ components of (5.10) are equivalent to our capacity constraints (5.7). For $j = b+1, \ldots, n$ we have $x_j^* = 0$ (i.e., activity $j$ is null in our nominal processing plan), and so components $r+1, \ldots, p$ of (5.10) simply express the fact that *cumulative activity levels for nonbasic activities must be nondecreasing.* Of course, the same is true for basic activities, and this can be expressed in terms of the centered processes $V_j(\cdot)$ as follows:

$$V_j(t) - V_j(0) \leq (t-s) x_j^* \qquad \text{for } j = 1, \ldots, b \text{ and } 0 \leq s \leq t, \tag{5.11}$$

which is identical to constraint (3.10) of the fluid model.

To form a Brownian approximation for our perfectly balanced stochastic processing network (Assumption 1 remains in force but will be relaxed shortly), we consider a small parameter $\varepsilon > 0$ and an arbitrary policy $T$, assuming that $\varepsilon$ and the initial inventory vector $Q(0)$ are such that

(5.12) all components of the $m$-vector $q = \varepsilon Q(0)$ are moderate in value,

then define the following scaled processes associated with $T$:

$$Y(t) = \varepsilon V(t/\varepsilon^2), \qquad t \geq 0, \tag{5.13}$$

$$Z(t) = \varepsilon Q(t/\varepsilon^2), \qquad t \geq 0, \tag{5.14}$$

$$U(t) = KY(t) = \varepsilon I(t/\varepsilon^2), \qquad t \geq 0. \tag{5.15}$$

Relating $\varepsilon$ to the scaling parameter $n$ in [6] via $\varepsilon = 1/\sqrt{n}$, one sees that our processes $Y$, $Z$ and $U$ are precisely analogous to the corresponding processes in [6], with one notable exception: here the scaled "unused capacity process" $U$ has been augmented to include another $n - b$ components representing cumulative activity levels for nonbasic activities; the treatment of multiclass queueing networks in [6] was framed in such a way that *all* processing activities were basic (i.e., all activities were conducted at positive levels in the nominal plan), and so the need for this crucial augmentation did not arise. Arguing exactly as in Section 5 of [6], one arrives at an approximating Brownian network whose key system equation is

$$Z(t) = q + X(t) + RY(t), \qquad t \geq 0, \tag{5.16}$$

where $X$ is the $m$-dimensional Brownian motion that approximates our scaled nominal workload process (5.4); that is, $X$ has zero drift and covariance matrix $\Sigma$ defined in terms of original model data via (5.5). Constraints (5.6) and (5.10) of the stochastic network model are reexpressed in the Brownian approximation as

$$Z(t) \geq 0 \qquad \text{for all } t \geq 0 \tag{5.17}$$

and

$$U(\cdot) \text{ is nondecreasing with } U(0) \geq 0. \tag{5.18}$$

The one aspect of (5.17) and (5.18) that deserves comment is the allowance of controls $Y$ such that $U(0) = KY(0)$ has one or more *strictly* positive components. This is appropriate, as in the approximation of multiclass queueing networks [6], because our diffusion scaling in (5.13)–(5.15) allows arbitrarily rapid increases in $U(\cdot)$ as $\varepsilon$ becomes small. In like fashion, the constraints (5.11), which express the requirement that *basic* activity levels be nonnegative, simply do not appear in the Brownian approximation, because those constraints become inconsequentially weak as $\varepsilon \searrow 0$ under diffusion scaling; this too is an observation made earlier in the treatment of queueing networks [6]. A final and obviously important constraint of our Brownian network approximation is that the chosen control $Y$ be adapted to the Brownian motion $X$, which expresses the requirement that activity levels up to time $t$ depend only on information available at $t$. This requirement, that flow management policies be approximately *nonanticipating*, has not even been mentioned in the exposition thus far, and a rigorous justification of its representation in the Brownian system model will not be attempted here.

Thus far attention has been restricted to a stochastic processing network whose data satisfy Assumption 1 exactly (i.e., a perfectly balanced network) and after rescaling time and state space by means of a small parameter $\varepsilon > 0$, we have arrived at an approximating Brownian network of the form (1.1)–(1.5), with $\theta = 0$, $\Sigma$ defined by (5.5) and the matrices $R$ and $K$ determined from data of the original model as in Sections 2 and 3. The situation where Assumption 1 holds *approximately* can be formalized as follows. Let us suppose that there exists a vector $\lambda^*$ such that Assumption 1 holds with $\lambda^*$ in place of $\lambda$ (again we denote by $x^*$ the optimal solution of the static planning problem, so $x^*$ satisfies $Rx^* = \lambda^*$ and $Ax^* = e$) and $\lambda$ is close to $\lambda^*$ in the following sense: there exists a small scalar $\varepsilon > 0$ such that (5.12) holds and moreover

$$\text{all components of the } m\text{-vector } \theta = (\lambda - \lambda^*)/\varepsilon \text{ are moderate in value.} \tag{5.19}$$

Again we call $x^*$ a vector of nominal activity levels, deviation controls are defined in terms of $x^*$ via (5.8) and the parameter $\varepsilon$ is used to rescale time and state space. As in Section 5 of [6], this leads to an approximating Brownian network identical to the one obtained earlier except that the Brownian motion $X$ has the drift vector $\theta$ defined in (5.19).

It is noteworthy that nonbasic activities, indexed by $j = b+1, \dots, n$, have no role in defining the covariance matrix $\Sigma$ in (5.5), because $x_j^* = 0$ for all such $j$. As noted earlier, a system manager must use those activities sparingly if excessive inventory build-ups are to be avoided, but nonbasic activities remain a relevant part of our approximating Brownian system model: their average input–output characteristics are still represented in the Brownian model as columns $b+1, \dots, n$ of the matrix $R$ in (5.18), and it will be shown by example in Section 6 that nonbasic activities may play an important role in dynamic control of the Brownian network.

As the final order of business in this section, let us consider the alternative scaling that gives a fluid limit or fluid approximation, as opposed to the diffusion or Brownian approximations considered up to this point. To simplify discussion, attention is restricted to a single system satisfying Assumption 1. Reasoning as in Sections 5 and 6 of [7] and in the papers cited there, one concludes that the desired scaling, using a small parameter $\varepsilon > 0$ such that (5.12) holds, is

$$\beta(t) = \varepsilon V(t/\varepsilon), \qquad t \geq 0, \tag{5.20}$$

$$z(t) = \varepsilon Q(t/\varepsilon), \qquad t \geq 0, \tag{5.21}$$

$$U(t) = K\beta(t) = \varepsilon I(t/\varepsilon), \qquad t \geq 0, \tag{5.22}$$

which leads to precisely the fluid model described in Section 3. The spatial scaling factor in (5.20)–(5.22) is exactly the same as in (5.13)–(5.15), but to achieve the fluid limit we use a weaker scaling of time in (5.20)–(5.22), so that one unit of scaled time in the fluid model corresponds to a time span of length $\varepsilon$ in the Brownian model. This distinction illuminates two equivalent characterizations of state space collapse obtained earlier: two states $q$ and $q'$ are equivalent in the Brownian model, meaning that either one can be *instantaneously* exchanged for the other, if and only if they communicate in the fluid model, meaning that either one can be *eventually* exchanged for the other.

**6. Examples and additional properties.** The workload process $W(t)$ in our reduced Brownian network is defined as $W(t) = MZ(t)$, and given the everyday meaning of the word "workload," one naturally expects all elements of $M$ to be nonnegative. This need not be true, however, as the following example shows. Consider a processing network with $r = 1$ and $m = n = 2$ (i.e., one resource or server, two materials and two activities). The first-order data are

$$\lambda = \begin{pmatrix} 3/2 \\ 1/2 \end{pmatrix}, \qquad R = \begin{bmatrix} 2 & 1 \\ 2 & -1 \end{bmatrix} \quad \text{and} \quad A = [1 \quad 1].$$

The key feature here is that activity 2 consumes one material and produces the other, whereas activity 1 consumes both materials at a relatively high rate; the two activities are equally expensive in terms of server capacity. These problem data satisfy Assumption 1, and the optimal solution of our static planning problem (2.1) and (2.2) is $x^* = (\frac{1}{2}, \frac{1}{2})$. The dual problem (2.3) and (2.4) has a unique optimal solution with $\mu^* = (\frac{3}{4}, -\frac{1}{4})$ and $\pi^* = 1$. Thus the reduced Brownian network has dimension $d = 1$ and the (unique) canonical basis matrix is $M = \mu^*$. Any other basis matrix for the linear space $\mathcal{M}$ would differ from this canonical choice by a scale factor, of course, and so the matrix $M$ that is used to define the workload process $W(t)$ via (1.8) cannot be taken nonnegative in this example.

Of course, the optimal dual variable $\mu_i^*$ tells us the rate at which $p$ increases due to an increase in the exogenous input rate $\lambda_i$. In this example, an increase in $\lambda_2$ actually *decreases* utilization $\rho$ of the single server because it enables greater use of activity 1. The following monotonicity assumption, which is satisfied by conventional queueing network models and is natural in many contexts, rules out such phenomena. It is a condition on the two matrices $R$ and $A$ that characterize first-order processing capabilities.

ASSUMPTION 2. If $\lambda > \lambda'$, $x \geq 0$ and $Rx = \lambda$, then there exists $x' \geq 0$ such that $Rx' = \lambda'$ and $Ax' \leq Ax$.

To paraphrase, Assumption 2 says that if one input stream is uniformly smaller than another in terms of average input rates, then it can be processed with uniformly lower long-run utilization rates for all resources. It is easy to show that no optimal solution $(\mu, \pi)$ of the dual problem (2.3) and (2.4) can include a negative $\mu_i$ value in this case, and so we have the following.

PROPOSITION 5. *Suppose that the first-order data ($R$, $A$ and $\lambda$) of an open processing network satisfy Assumption 1, and that $R$ and $A$ further satisfy Assumption 2. Then the canonical basis matrix $M$ in (2.7) is non-negative.*

Recall that in (2.8) we established the notation $R = (H, J)$, where the columns in $H$ and $J$ correspond to basic and nonbasic activities, respectively. The following condition is satisfied by virtually all interesting models.

ASSUMPTION 3. $H$ has full row dimension $m$, and thus there exists a $b \times m$ matrix $H^+$ (a right inverse for $H$) satisfying $HH^+ = I$.

PROPOSITION 6. *Suppose that the first-order data ($R$, $A$ and $\lambda$) of an open processing network satisfy Assumptions 1 and 3. Then the canonical basis matrix $M$ in (2.7) can be decomposed as $M = \Pi B H^+$, where the non-negative matrices $\Pi$ and $B$ are defined in (2.7) and (2.8), respectively.*

PROOF. From (2.9) we have $MH = \Pi B$, and right-multiplying both sides of this equation by $H^+$ gives the desired representation. □

To grasp the contents of Assumption 3, it is perhaps most illuminating to examine a model for which it does *not* hold. Consider a system with $r = 1$, $m = 2$ and $n = 3$ (i.e., one server, two materials and three activities). Let the first-order data be

$$\lambda = \begin{pmatrix} 1 \\ 1 \end{pmatrix}, \qquad R = \begin{bmatrix} 1 & 4/3 & 0 \\ 1 & 0 & 4/3 \end{bmatrix} \quad \text{and} \quad A = [1 \quad 1 \quad 1].$$

These data satisfy Assumption 1, the unique optimal solution of (2.1) and (2.2) being $x^* = (1, 0, 0)$. The dual problem (2.3) and (2.4) has exactly two basic feasible solutions, and thus there are two generalized cut constraints: one of the basic feasible solutions is $\mu^1 = (\frac{3}{4}, \frac{1}{4})$ and $\pi^1 = 1$; the other is $\mu^2 = (\frac{1}{4}, \frac{3}{4})$ and $\pi^2 = 1$. With the vector $\lambda$ specified above, both of those basic feasible solutions happen to be optimal, or equivalently, both of the generalized cut constraints are binding for the given $\lambda$. Thus we have a "reduced Brownian network" of dimension $d = 2$, which shows that *the effective system dimension $d$ for a Brownian network may be strictly larger than the number $r$ of servers or resources.* Also, in this example, $H$ is a $2 \times 1$ matrix consisting of just the first column of $R$, so Assumption 3 does not hold. The salient characteristic of this example is a sort of degeneracy that allows two exogenous input materials to be processed in exactly the required relative volumes using a single activity.

This example provides a good opportunity to examine the role of nonbasic activities in a Brownian system model. For concreteness, assume that the system manager wants to minimize total inventory, represented by $Z_1(t) + Z_2(t)$. One ultimately finds that there exists a *pathwise-optimal* control policy for the Brownian model, which means that the policy minimizes $Z_1(t) + Z_2(t)$ for all $t \geq 0$ with probability 1. Because the equivalent workload formulation has dimension $d = m = 2$, we can simply retain the two-vector $Z(t)$ as our state descriptor and consider deviation controls $Y = (Y_1, Y_2, Y_3)$ that might be employed. An increase in $Y_1$ corresponds to decreasing activity 1 from its nominal rate of $x_1^* = 1$ (which fully absorbs all server capacity), whereas *decreases* in $Y_2$ and $Y_3$ correspond to insertion of activities 2 and 3, respectively, both of which are null in the nominal plan. The policy constraints are that $Y_2$ and $Y_3$ be nonincreasing and that $U = AY = Y_1 + Y_2 + Y_3$ be nondecreasing. Of course, $Y$ must also be chosen so that the inventory vector $Z$ remains nonnegative. Looking at the second column of $R$, we see that a decrease of $\delta$ in $Y_2$ (i.e., inserting $\delta$ units of activity 2), decreases $Z_1$ by $4\delta/3$ relative to where it would have been under the nominal processing plan, but this must be accompanied by an increase of $\delta$ in $Y_1$ (i.e., deleting $\delta$ units of activity 1 relative to the nominal plan) in order to keep $U$ nondecreasing. The latter action increases both $Z_1$ and $Z_2$ by $\delta$ relative to the values they would have had under the nominal plan, and so the net effect is to displace $Z = (Z_1, Z_2)$ by $(\delta - 4\delta/3, \delta)$ and hence increase $Z_1 + Z_2$ by $2\delta/3$. To minimize $Z_1 + Z_2$ in the Brownian model, one then finds that the following policy is optimal: exert no control (i.e., leave the three-vector $Y$ at its current value) when both components of $Z$ are strictly positive; increase $Y_1$ and decrease $Y_2$ by equal amounts when the boundary $Z_2 = 0$ is struck, using the minimal control quantities required to keep $Z_2 \geq 0$; and apply the controls $Y_1$ and $Y_3$ in like fashion at the boundary $Z_1 = 0$ so as to keep $Z_1 \geq 0$. In terms of our original processing network, this obviously means that the system manager should use all server capacity on activity 1 when both of the required inputs are available (i.e., when neither buffer is empty), but substitute activity 2 for activity 1 when buffer 2 is empty, and substitute activity 3 for activity 1 when buffer 1 is empty. Under heavy traffic conditions this policy actually uses very little server capacity in either activity 2 or activity 3, and yet the availability of these substitute nonbasic activities (as alternatives to simply idling the server) has a substantial effect on system behavior under diffusion scaling.

Suppose that $R$ and $A$ are the same as in the last example but $\lambda = (\frac{1}{4}, \frac{5}{4})$. Again Assumption 1 is satisfied, the unique solution of (2.1) and (2.2) being $x^* = (\frac{1}{4}, 0, \frac{3}{4})$, and now the *unique* solution of the dual problem (2.3) and (2.4) is $(\mu^2, \pi^2)$. That is, only the second of our two generalized cut constraints is binding with the alternative $\lambda$ vector, and the effective system dimension is $d = 1$. Now $H$ is the $2 \times 2$ nonsingular matrix consisting of the first and third columns of $R$, so the unique right inverse is

$$H^+ = H^{-1} = \begin{bmatrix} 1 & 0 \\ -3/4 & 3/4 \end{bmatrix}.$$

Also, $\Pi = (\pi^2) = (1)$ and $B = (1, 1)$, so Proposition 6 gives $M = \Pi B H^+ = (\frac{1}{4}, \frac{3}{4}) = \mu^1$.

For a general interpretation of the decomposition in Proposition 6, it is useful to consider an inventory vector $z \in \mathbb{R}^m_+$ and define the corresponding workload vector

$$w = Mz = \Pi B H^+ z. \tag{6.1}$$

Suppose that the right-inverse $H^+$ is nonnegative. Noting that the vector $x = H^+ z$ satisfies $Hx = z$, we interpret $H^+ z$ as a vector of basic activity levels that can be used to process the inventories in $z$ to completion, and $BH^+ z$ as a vector of capacities consumed by that program of activities, with one component for each server or resource. That is, the $r$ components of $BH^+ z$ show the total amounts of work required from the $r$ different servers (expressed in time units) to process an initial inventory of $z$. But the minimal description $w$ of system workload in (6.1) may actually be coarser than that because of resource pooling: each row of the *pooling matrix* $\Pi$ specifies a positive linear combination of the $r$ different server workloads; we interpret the corresponding component of $w$ as the workload for a resource pool whose members are the servers having strictly positive coefficients. Laws [14, 15] emphasized the resource pooling that occurs in queueing networks where alternate routing capabilities exist, and one of his examples will be reviewed in Section 7 below; another example of resource pooling due to alternate routing is discussed in Section 5 of [10]. In a general open processing network there is no guarantee that the right-inverse $H^+$ can be chosen nonnegative, but one can still interpret components of the vector $H^+ z$ in (6.1) as *changes* in basic activity levels over a span of time, relative to the nominal processing plan, that are required to decrease the system's inventory vector by $z$.

In a general processing network, the term "alternate routing" may not be applicable, but one can still say that resource pooling occurs when the set of basic activities is rich enough to provide alternative means of processing the exogenous input flows. The following example, which is closely related to the one discussed immediately above, illustrates this phenomenon. Consider a two-server system with two materials, three activities and first-order data

$$\lambda = \begin{pmatrix} 5/4 \\ 19/12 \end{pmatrix}, \qquad R = \begin{bmatrix} 1 & 4/3 & 0 \\ 1 & 0 & 4/3 \end{bmatrix} \quad \text{and} \quad A = \begin{bmatrix} 1 & 1 & 0 \\ 0 & 0 & 1 \end{bmatrix}.$$

These satisfy Assumption 1, the unique solution of our static planning problem (2.1) and (2.2) being $x^* = (\frac{1}{4}, \frac{3}{4}, 1)$. The dual problem (2.3) and (2.4) also has a unique solution, which is $\mu^* = (\frac{9}{16}, \frac{3}{16})$ and $\pi^* = (\frac{3}{4}, \frac{1}{4})$. Because all three activities are basic, we have $H = R$ and $B = A$. One feasible choice of the right-inverse $H^+$ is

$$H^+ = \begin{bmatrix} 0 & 1 \\ 3/4 & 0 \\ 0 & 3/4 \end{bmatrix},$$

and in the end we have

$$M = \Pi B H^+ = \pi^* A H^+ = (9/16, 3/16) = \mu^*.$$

The pooling matrix $\Pi = \pi^* = (\frac{3}{4}, \frac{1}{4})$ shows that the two servers function as a single capacity pool for purposes of our Brownian system model, with server 1 contributing three times as much "effective capacity" as server 2.

**7. An example due to Laws.** Consider the processing network pictured in Figure 1, where six servers (represented by circles) are arranged in a $2 \times 3$ array. The open-ended rectangles represent buffers in which jobs of nine different classes are stored. There are external arrivals into buffers 1 and 2 only (at rates $\lambda_1$ and $\lambda_2$, resp.), and the system manager has discretion as to how new arrivals will be processed: class 1 jobs can be processed by either server 1 or server 4, class 2 jobs can be processed by any of servers 1, 2 or 3, and as the arrows in Figure 1 show, the routing of a job after its first processing operation is completely determined by which server is chosen to perform that initial operation. Thus there are a total of twelve processing activities available to the system manger: two ways to process class 1 jobs, three ways to process class 2 jobs and one way to process each of the other seven job classes. Assuming that each of the twelve service time distributions has mean 1, the input–output matrix $R$ and resource consumption matrix $A$ for this network are as follows:

$$R=\begin{bmatrix}
1 & 1 & & & & & & & & & & \\
 & & 1 & 1 & 1 & & & & & & & \\
-1 & & & & & 1 & & & & & & \\
 & & & & & -1 & 1 & & & & & \\
 & & -1 & & & & & 1 & & & & \\
 & & & -1 & & & & & 1 & & & \\
 & & & & -1 & & & & & 1 & & \\
 & -1 & & & & & & & & & 1 & \\
 & & & & & & & & & & -1 & 1
\end{bmatrix}$$

and

$$A=\begin{bmatrix}
1 & & 1 & & & & & & & & & \\
 & & & 1 & & 1 & & & & & & \\
 & & & & 1 & & 1 & & & & & \\
 & 1 & & & & & & 1 & & & & \\
 & & & & & & & & 1 & & 1 & \\
 & & & & & & & & & 1 & & 1
\end{bmatrix}.$$

In this representation, activities 1 and 2 correspond to the processing of class 1 jobs by servers 1 and 4, respectively, activities 3 through 5 correspond to the processing of class 2 jobs by servers 1 through 3, respectively, and activities 6 through 12 correspond to the processing of classes 3 through 9, respectively, by the one server which is qualified to process that class. (E.g., only server 6 can process class 9, and that is activity 12 in our system representation.)

This $2\times 3$ example was originally introduced by Laws [14, 15] and further discussed in the survey paper by Kelly and Laws [13], but Laws' representation of the system differs from the one presented here. Laws assumes that routing decisions must be made immediately when new jobs arrive, whereas we have assumed that new arrivals of either class 1 or class 2 can be held in a common buffer, with routing decisions delayed until the system manager commits a server to their initial processing operation. This distinction causes Laws to define three more buffers than shown in Figure 1 (two buffers rather than one for the horizontal input stream and three buffers rather than one for the vertical input stream) but the distinction between immediate and delayed routing is inconsequential in heavy traffic, and the reduced Brownian system model to be derived here (i.e., our equivalent workload formulation) is identical to the one derived by Laws.

A more consequential difference is that Laws formulates his analog of the primal linear program (2.1) and (2.2) in terms of complete routes rather than elementary processing activities. For example, in Figure 1 there are exactly two complete routes available to class 1 jobs and three complete routes available to class 2 jobs, so the coefficient matrices appearing in Laws' analog of (2.1) and (2.2) have only five columns. For current purposes it is not necessary to review his LP formulation in detail, but the following remarks are called for. First, Laws' theory is restricted to processing networks where each activity involves a single server processing a single job class (i.e., no activity may involve simultaneous use of two or more servers, nor may it involve multiple inputs), with no uncertainty as to where the job can or must go after that processing is complete, and it is not obvious how to extend his style of LP formulation to the more general class of networks considered here. Second, our formulation (2.1) and (2.2) involves exactly the same first-order data ($R$, $A$, $\lambda$) that appear in the natural fluid approximation for the processing network (see Section 3), and that parallelism has a number of advantages. Despite these differences, the two LP formulations are essentially equivalent for the class of models considered by Laws, and we shall cite various of his findings in the paragraphs that follow, translating as necessary into the framework and notation developed in this paper.

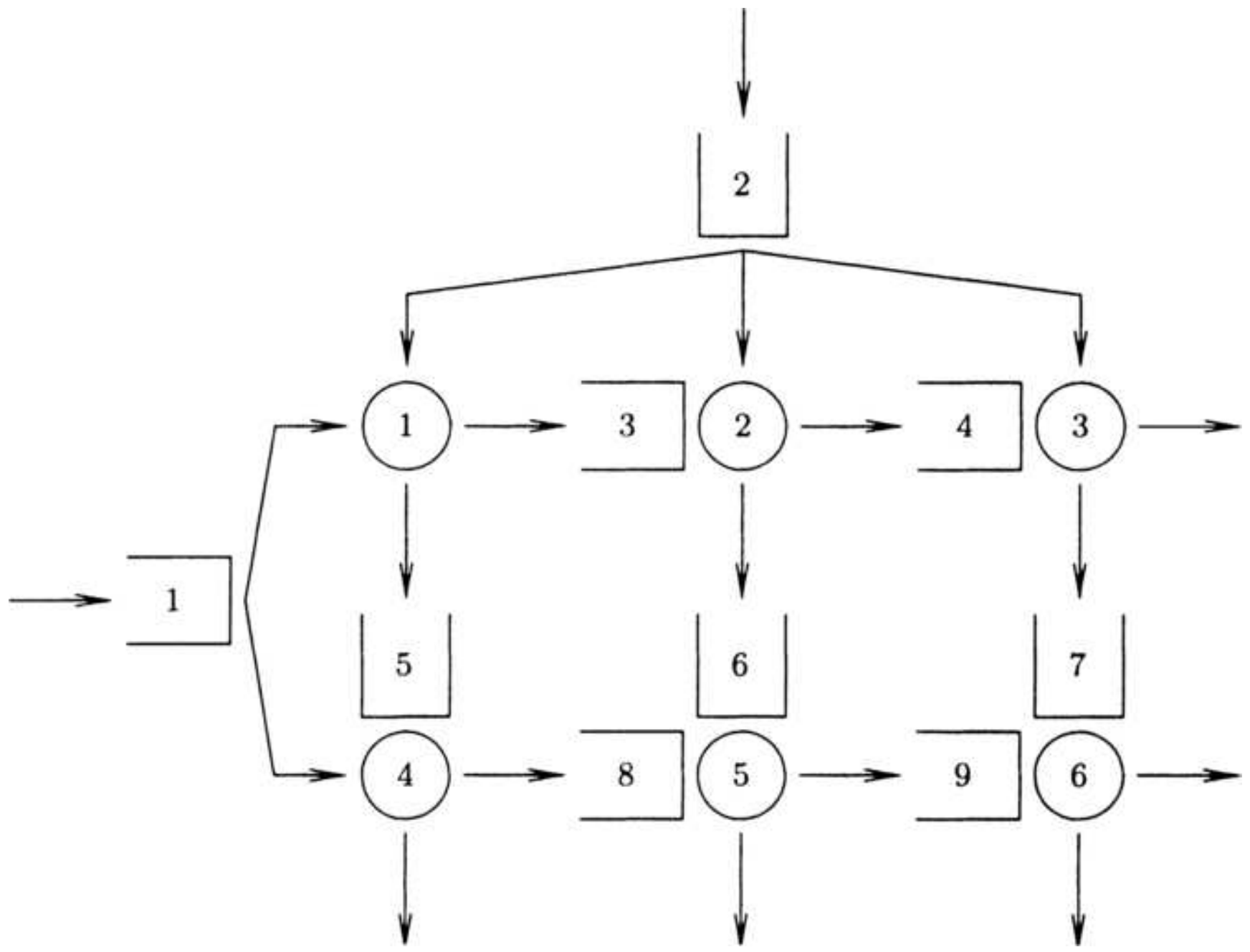

FIG. 1. *The* $2 \times 3$ *network.*

Assuming hereafter that the external arrival rates satisfy $\frac{1}{2}\lambda_1 + \frac{1}{3}\lambda_2 = 1$, we find that the static planning problem (2.1) and (2.2) has a unique solution, namely: $\rho^* = 1$; $x_j^* = \frac{1}{2}\lambda_1$ for $j = 1, 2, 6, 7, 11$ and 12 (these are activities which process what were originally class 1 arrivals); and $x_j^* = \frac{1}{3}\lambda_2$ for $j = 3, 4, 5, 8, 9$ and 10 (these are activities which process what were originally class 2 arrivals). Thus all twelve activities are basic, which implies that the set $\mathcal{D}^*$ of optimal solutions for our dual linear program (2.3) and (2.4) consists precisely of pairs $(\mu, \pi)$ satisfying

$$\mu R = \pi A, \qquad \pi e = 1 \quad \text{and} \quad \pi \geq 0. \tag{7.1}$$

Thus, using the language of "workload contributions" and "relative server capacities" that was proposed in Section 2, each optimal solution of the dual problem corresponds to a set of nonnegative relative server capacities $\pi_k$ that sum to one and have the following property: for each buffer $i$, the sum of the server capacities consumed in processing to completion a unit of class $i$ material is equal to a common value $\mu_i$, regardless of which route (i.e., regardless of which available processing sequence) may be used to process the material. In similar fashion, Laws [14] interprets solutions $(\mu, \pi)$ of (7.1) in terms of a tolling scheme: $\pi_k$ represents the fee charged for server $k$ capacity, and $\mu_i$ is then the cost incurred (i.e., the sum of all fees paid) in processing to completion a class $i$ job; the key requirement is that the sum of the fees paid must be the same over all routes available to any given job class.

A detailed analysis shows that the set $\mathcal{D}^*$ defined by (7.1) has six distinct extreme points, which are precisely the alternative basic optimal solutions $(\mu^l, \pi^l)$ identified in Section 2, as follows:

$$\begin{aligned}
6\mu^1 &= (3,2,1,1,0,2,1,3,1) \quad \text{and} \quad 6\pi^1 = (2,0,1,0,2,1);\\
6\mu^2 &= (3,2,1,0,0,1,2,3,2) \quad \text{and} \quad 6\pi^2 = (2,1,0,0,1,2);\\
6\mu^3 &= (3,2,3,1,2,0,1,1,1) \quad \text{and} \quad 6\pi^3 = (0,2,1,2,0,1);\\
6\mu^4 &= (3,2,2,0,1,0,2,2,2) \quad \text{and} \quad 6\pi^4 = (1,2,0,1,0,2);\\
6\mu^5 &= (3,2,3,2,2,1,0,1,0) \quad \text{and} \quad 6\pi^5 = (0,1,2,2,1,0);\\
6\mu^6 &= (3,2,2,2,1,2,0,2,0) \quad \text{and} \quad 6\pi^6 = (1,0,2,1,2,0).
\end{aligned}$$

In this example, because all activities are basic, the matrices $H$ and $B$ introduced in (2.8) are identical with $R$ and $A$, respectively. It is easy to verify that $R$ has full row dimension, so Assumption 3 is satisfied, and hence one has $\mu^l = \pi^l A R^+$ for each $l = 1, \ldots, 6$ (see Proposition 6). Because the $6 \times 9$ matrix $AR^+$ has full row dimension, any given subset of the vectors $\mu^1, \ldots, \mu^6$ has the same dimension as the corresponding subset of the vectors $\pi^1, \ldots, \pi^6$. It can be verified that $\pi^1, \pi^2$ and $\pi^3$ are linearly independent and that their span includes $\pi^4, \pi^5$ and $\pi^6$. Thus the dimension of our equivalent workload formulation is $d = 3$; our "canonical choice" of basis matrix is

$$M = \begin{bmatrix} \mu^1 \\ \mu^2 \\ \mu^3 \end{bmatrix} = \frac{1}{6} \begin{bmatrix} 3 & 2 & 1 & 1 & 0 & 2 & 1 & 3 & 1 \\ 3 & 2 & 1 & 0 & 0 & 1 & 2 & 3 & 2 \\ 3 & 2 & 3 & 1 & 2 & 0 & 1 & 1 & 1 \end{bmatrix}, \tag{7.2}$$

which has the representation

$$(7.3) \qquad M = \Pi A R^{+} \qquad \text{where } \Pi = \begin{bmatrix} \pi^1 \\ \pi^2 \\ \pi^3 \end{bmatrix} = \frac{1}{6}\begin{bmatrix} 2 & 0 & 1 & 0 & 2 & 1 \\ 2 & 1 & 0 & 0 & 1 & 2 \\ 0 & 2 & 1 & 2 & 0 & 1 \end{bmatrix}.$$

To repeat, we have found that $d < L^*$ in this example, and thus our "canonical choice" of the basis matrix $M$ is nonunique, depending on the particular order in which the extreme points $(\mu^l, \pi^l)$ of $\mathcal{D}$ were enumerated.

In Section 6 of [13], Kelly and Laws state that a suitable definition of workload for the Brownian model of this network is a three-dimensional process $W(t) = MZ(t)$, where $M$ is computed from precisely the same pooling matrix $\Pi$ displayed in (7.3). Their matrix $M$ is actually different from (7.2) because of the difference in system representation referred to earlier, but their three-dimensional workload process is entirely equivalent to the one derived here. The justification for that equivalent workload formulation is contained in Laws' Ph.D. dissertation [14], where virtually all the results presented in this paper were developed for the restricted class of networks that he considered (see above). In fact, Laws' treatment of those networks went beyond what is given in this paper, because he characterized precisely the "nonbottle-neck servers" whose influence is negligible in the heavy traffic limit, whereas such servers have been deleted with only a passing remark in our treatment. Apart from the examples and exposition contained in the Kelly-Laws survey paper [13], the only part of Laws' general theory that has been published in the open literature is that concerned with networks whose Brownian models have a one-dimensional equivalent workload formulation [15]. (In our notation, this is the case $d = 1$.) That is, Laws' paper [15] deals with networks of the restricted type described earlier in this section where, moreover, just one "generalized cut constraint" is binding in the heavy traffic limit; one must look in his Ph.D. dissertation [14] to find his deep and illuminating derivation of equivalent workload formulations when several such constraints are binding in the heavy traffic limit.

Turning to the general interpretation of the optimal dual variables $\mu$ and $\pi$ that was developed in Section 3, one sees that each of the vectors $\pi^1, \ldots, \pi^6$ displayed above corresponds to a different pair of servers that may be rendered "noncritical" by a perturbation of the following sort. Starting with large amounts of material in buffers 1 and 2, but all other buffers empty, the processing plan which minimizes time-to-emptiness $\tau$ (with exogenous inputs turned off) is one that keeps all servers busy throughout the interval $[0, \tau]$. Now imagine that a small amount $\delta_i$ of additional material is placed in each buffer $i$, assuming for simplicity that the entire vector $\delta = (\delta_i)$ is nonnegative. To see which servers remain "critical" in minimizing time-to-emptiness, we compute $\mu^l \delta$ for each $l = 1, \ldots, 6$ and focus on the value of $l$ for which that product is maximal. For concreteness, suppose it is $l = 1$.

Then the processing plan which minimizes time-to-emptiness is one in which both server 2 and server 4 experience some idleness; one knows this because $\pi_2^1 = \pi_4^1 = 0$. Thus added capacity from those servers would have no impact on MTTE, but the four positive values $\pi_k^1$ show the effect (in terms of reducing MTTE) of added capacity for the servers who remain critical after the perturbation. Readers are invited to hypothesize specific perturbations $\delta$ to see whether they can foresee which servers are rendered noncritical without actually computing $\mu'\delta$ for all $1 = 1, \dots, 6$.

**Acknowledgments.** I am indebted to Neil Laws for catching a major error in the original version of this paper and to Frank Kelly and Sunil Kumar for their helpful comments.

Graduate School of Business
Stanford University
Stanford, California 94305-5015
USA
E-mail: harrison_michael@gsb.stanford.edu